\documentclass[12pt]{article} 
 \usepackage[margin=1in]{geometry}
\usepackage{setspace}
\usepackage{graphicx}	
\usepackage{subcaption}
\usepackage{float} 
\usepackage[toc]{appendix}
\usepackage[draft]{todonotes} 
\usepackage{enumerate}	
\usepackage{paralist}	
\usepackage{amsmath, amsthm, amssymb,mathtools}
\usepackage{thm-restate}
\usepackage{wasysym}
\usepackage[pdftex,hypertexnames=false,linktocpage=true,hidelinks]{hyperref}
\usepackage{cleveref}
\usepackage{darkmode}
\usepackage{multirow}
\usepackage{ circuitikz }

\newtheorem{thm}{Theorem}[section]

\newtheorem{lemma}[thm]{Lemma}

\theoremstyle{definition}

\newtheorem{defn}[thm]{Definition}
\theoremstyle{definition}

\newtheorem{obs}[thm]{Observation}

\renewcommand\emptyset{\varnothing}

\def\MIS{\mathrm{MIS}}
\def\NIMIS{\mathrm{NIMIS}}

\newcommand{\Grid}[2]{G_{#1 \times #2}}
\newcommand{\Band}[2]{FC_{#1 \times #2}}
\newcommand{\Tube}[2]{TC_{#1 \times #2}}

\newcommand{\Mobius}[2]{M_{#1 \times #2}}

\newcommand\Stab[1]{Stab(#1)}
\newcommand{\Mod}[1]{\ (\mathrm{mod}\ #1)}

\usepackage{tikz}
\usepackage{xcolor}

\usepackage{bbm}

\definecolor{amber}{rgb}{1.0, 0.75, 0.0}

\definecolor{forest}{rgb}{0.0, 0.5, 0.0}

\definecolor{cadmium}{rgb}{0.93, 0.53, 0.18}

\definecolor{bittersweet}{rgb}{1.0, 0.44, 0.37}

\definecolor{byzantine}{rgb}{0.74, 0.2, 0.64}

\definecolor{brilliantrose}{rgb}{1.0, 0.33, 0.64}

\definecolor{caribbeangreen}{rgb}{0.0, 0.8, 0.6}

\definecolor{electriccyan}{rgb}{0.0, 1.0, 1.0}

\definecolor{periwinkle}{rgb}{0.8, 0.8, 1.0}

\definecolor{steelblue}{rgb}{0.27, 0.51, 0.71}

\title{Statistics of maximal independent sets in grid-like graphs}
  \author{Levi Axelrod\thanks{Department of Mathematics and Statistics, Binghamton University, Binghamton, NY, E-mail: {\tt laxelro1@binghamton.edu}.} \and Nathan Bickel\thanks{Department of Mathematics, Iowa State University, Ames, IA, E-mail: {\tt nbickel@iastate.edu}. Research of this author was supported by NSF Grant DMS-1950583.} \and Anastasia Halfpap\thanks{Department of Mathematics, Iowa State University, Ames, IA, E-mail: {\tt ahalfpap@iastate.edu}. Research of this author was supported in part by NSF FRG DMS-2152490.} \and Luke Hawranick\thanks{School of Mathematical and Data Sciences, West Virginia University, Morgantown, WV, E-mail: {\tt lh00022@mix.wvu.edu}. Research of this author was supported by NSF Grant DMS-1950583.} \and Alex Parker\thanks{Department of Mathematics, Iowa State University, Ames, IA, E-mail: {\tt abparker@iastate.edu}.} \and Cole Swain\thanks{Department of Biology, Chemistry, Mathematics and Computer Science, University of Montevallo, Montevallo, AL, E-mail: {\tt cxswain@go.olemiss.edu}. Research of this author was supported by NSF Grant DMS-1950583.}}

\begin{document}

\maketitle

\begin{abstract}
An independent set $I$ in a graph $G$ is \textit{maximal} if $I$ is not properly contained in any other independent set of $G$.  The study of maximal independent sets (MIS's) in various graphs is well-established, often focusing upon enumeration of the set of MIS's. For an arbitrary graph $G$, it is typically quite difficult to understand the number and structure of MIS's in $G$; however, when $G$ has regular structure, the problem may be more tractable. One class of graphs for which enumeration of MIS's is fairly well-understood is the rectangular grid graphs $G_{m\times n}$.

We say a graph is \textit{grid-like} if it is locally isomorphic to a square grid, though the global structure of such a graph might resemble a surface such as a torus or M\"obius strip. We study the properties of MIS's in various types of grid-like graphs, in particular determining parity of the set of MIS's, average size of MIS's, and number of pairwise non-isomorphic MIS's in various grid-like graphs.
\end{abstract}
\section{Introduction}\label{intro}

Given a graph $G$, a set $I \subseteq V(G)$ is called an $\textit{independent set}$ if no two vertices in $I$ are adjacent. A \textit{maximal independent set} (\textit{MIS}) $I$ in a graph $G$ is an independent set that is not a proper subset of another independent set. 
We say $I \subseteq V(G)$ \textit{covers} a vertex $u$ if $I$ contains either $u$ or a neighbor of $u$. Thus, $I$ is an MIS only if it covers every vertex in the graph. 
For a graph $G$, we denote by $\mathrm{MIS}(G)$ the collection of maximal independent sets in $G$.

Maximal independent sets are a structure of well-established interest, with many authors focusing on enumeration of MIS's in different contexts. For an arbitrary graph $G$, it is generally difficult to exactly describe the elements of $\MIS(G)$; determining $|\mathrm{MIS}(G)|$ is a $\#\mathcal{P}$-complete problem, even if $G$ is bipartite \cite{Vadhan}. However, the understanding of $\MIS(G)$ can become much more approachable if $G$ has some regular structure. One setting in which MIS's are fairly well-understood is the rectangular \textit{grid graphs}. 

\def\rows{4}
\def\cols{6}

\begin{defn}
    The $m \times n$ \textit{grid graph} $\Grid{m}{n}$ is the graph with $$V(\Grid{m}{n}):=\{(i,j) : 1 \leq j \leq m, 1 \leq i \leq n\}$$ and $$E(\Grid{m}{n}):=\{ \left\{ \left(i, j\right), 
    \left(i', j'\right) \right\} : (i, j), (i',j') \in V(\Grid{m}{n}), |i - i'| + |j - j'| = 1\}.$$ More concisely, $\Grid{m}{n}$ is the Cartesian product of $P_m$ and $P_n$, denoted $P_m \Box P_n$. (Throughout the paper, we use $P_n$ to denote the path on $n$ \textit{vertices}.)
\end{defn}

    \begin{figure}[h!]
    \begin{center}
    \begin{tikzpicture}
    
        \foreach \row in {1,...,\rows} {
            \foreach \col in {1,...,\cols} {
                \node[draw, circle, inner sep=3pt] (N\row\col) at (\col, -\row) {};
            }
        }
        
        \foreach \row in {1,...,\rows} {
            \foreach \col in {1,...,\numexpr\cols-1\relax} {
                \pgfmathtruncatemacro{\nextcol}{\col + 1}
                \draw (N\row\col) -- (N\row\nextcol);
            }
        }
        
        \foreach \col in {1,...,\cols} {
            \foreach \row in {1,...,\numexpr\rows-1\relax} {
                \pgfmathtruncatemacro{\nextrow}{\row + 1}
                \draw (N\row\col) -- (N\nextrow\col);
            }
        }

    \end{tikzpicture}
    \end{center}
    \caption{The grid graph $G_{4 \times 6}$}~\label{grid fig}
    \end{figure}
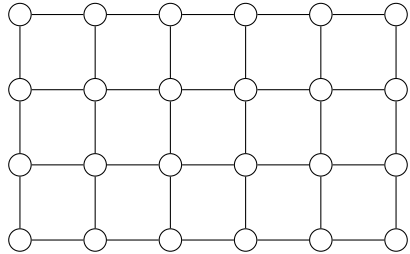

Study of MIS's in grids is a problem of established interest, both for pure mathematical and applied reasons. Various authors have studied the enumeration problem for rectangular grids, with Euler~\cite{Euler2} and Ortiz and Villanueva~\cite{Ortiz} giving explicit formulas in some cases. Oh~\cite{BigOh} gave a general formula for $|\MIS(\Grid{m}{n})|$ using the state matrix recursion algorithm. For related enumeration questions on grid graphs, see for instance~\cite{Transfer, OptimalDominating}.

The goal of this paper is to study the structure of maximal independent sets in grids and a variety of graphs with similar structure to grids. Collectively, we call such graphs \textit{grid-like}.

\begin{defn}\label{gridlike definition}
    Let $V_{m \times n} = \{(i,j) : 1 \leq j \leq m, 1 \leq i \leq n\}$ and 
    \[E_{m \times n} = \{ \left\{ \left(i, j\right), 
    \left(i', j'\right) \right\} : (i, j), (i', j') \in V(\Grid{m}{n}), |i - i'| + |j - j'| = 1\}.\] 
    The $m \times n$ grid-like graphs we study are:
    \begin{enumerate}
        \item The \textit{fat grid cylinder} $\Band{m}{n}$, with vertex set $V_{m \times n}$ and edge set 
        \[E_{m\times n} \cup \{\{(1,j), (n,j)\}: 1 \leq j \leq m\}\]
        \item The \textit{thin grid cylinder} $\Tube{m}{n}$, with vertex set $V_{m \times n}$ and edge set
        \[E_{m\times n} \cup \{\{(i,1), (i,m)\}: 1 \leq i \leq n\}\]
        \item The \textit{grid M\"obius strip} $\Mobius{m}{n}$, with vertex set $V_{m \times n}$ and edge set 
        \[E_{m \times n} \cup \{\left\{ \left(1, i\right), \left(n, m+1-i\right) \right\} : 1 \le i \le m\}.\]
    \end{enumerate}
\end{defn}

    We note that just as $\Grid{m}{n} = P_m \Box P_n$, we can more concisely describe $\Band{m}{n} = P_m \Box C_n$  and $\Tube{m}{n} = C_m \Box P_n$. Sometimes, this concise Cartesian product definition will suffice, but often we will wish to refer explicitly to the labeled vertices of grid-like graphs. As suggested by our naming conventions, the grid-like graphs which we study can be thought of visually as embeddings of the grid into various surfaces. However, this topological viewpoint will not be featured in our approach (beyond providing additional motivation for the study of these families and intuition for their symmetries). 

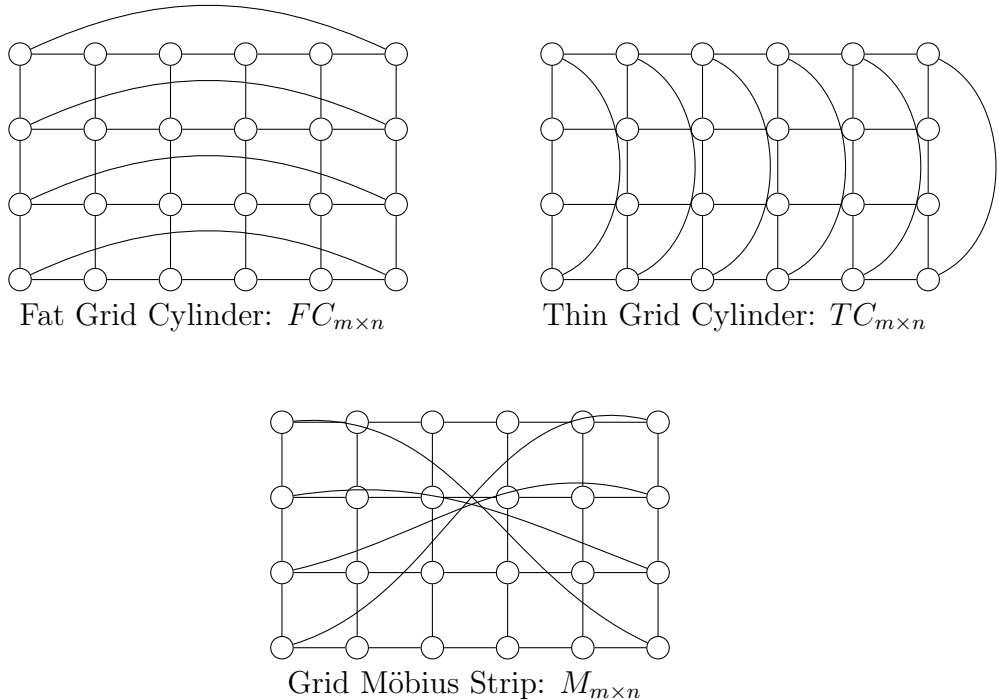
\begin{figure}[h!]

    \begin{center}
    
    \begin{minipage}[t]{0.45\linewidth}
        \begin{flushright}
        \begin{tikzpicture}
    
            \foreach \row in {1,...,\rows} {
                \foreach \col in {1,...,\cols} {
                    \node[draw, circle, inner sep=3pt] (N\row\col) at (\col, -\row) {};
                }
            }
            
            \foreach \row in {1,...,\rows} {
                \foreach \col in {1,...,\numexpr\cols-1\relax} {
                    \pgfmathtruncatemacro{\nextcol}{\col + 1}
                    \draw (N\row\col) -- (N\row\nextcol);
                }
                \draw (N\row1) to[out=25, in=155] (N\row\cols);
            }
            
            \foreach \col in {1,...,\cols} {
                \foreach \row in {1,...,\numexpr\rows-1\relax} {
                    \pgfmathtruncatemacro{\nextrow}{\row + 1}
                    \draw (N\row\col) -- (N\nextrow\col);
                }
            }
        
            \node at (\cols/1.75, -\rows-0.5) {Fat Grid Cylinder: $\Band{m}{n}$};
            
        \end{tikzpicture}
        \end{flushright}
    \end{minipage}
    \hfill
    \begin{minipage}[t]{0.45\linewidth}
        \begin{tikzpicture}
    
            \foreach \row in {1,...,\rows} {
                \foreach \col in {1,...,\cols} {
                    \node[draw, circle, inner sep=3pt] (N\row\col) at (\col, -\row) {};
                }
            }
            
            \foreach \row in {1,...,\rows} {
                \foreach \col in {1,...,\numexpr\cols-1\relax} {
                    \pgfmathtruncatemacro{\nextcol}{\col + 1}
                    \draw (N\row\col) -- (N\row\nextcol);
                }
            }
            
            \foreach \col in {1,...,\cols} {
                \foreach \row in {1,...,\numexpr\rows-1\relax} {
                    \pgfmathtruncatemacro{\nextrow}{\row + 1}
                    \draw (N\row\col) -- (N\nextrow\col);
                }
                \draw (N1\col) to[out=-25, in=25] (N\rows\col);
            }
        
            \node at (\cols/1.75, -\rows-0.5) {Thin Grid Cylinder: $\Tube{m}{n}$};
            
        \end{tikzpicture}
    \end{minipage}
    
    \vspace{0.5cm}

    \begin{minipage}[t]{0.45\linewidth}
    \begin{center}
        \begin{tikzpicture}
    
            \foreach \row in {1,...,\rows} {
                \foreach \col in {1,...,\cols} {
                    \node[draw, circle, inner sep=3pt] (N\row\col) at (\col, -\row) {};
                }
                
                \pgfmathtruncatemacro{\col}{\rows + 1 - \row}
                \pgfmathsetmacro{\angleout}{\row*\rows-\rows+5}
                \pgfmathsetmacro{\anglein}{150 + \angleout}
                \draw (N\row1) to[out=\angleout, in=\anglein] (N\col\cols);
            }
            
            \foreach \row in {1,...,\rows} {
                \foreach \col in {1,...,\numexpr\cols-1\relax} {
                    \pgfmathtruncatemacro{\nextcol}{\col + 1}
                    \draw (N\row\col) -- (N\row\nextcol);
                }
            }
            
            \foreach \col in {1,...,\cols} {
                \foreach \row in {1,...,\numexpr\rows-1\relax} {
                    \pgfmathtruncatemacro{\nextrow}{\row + 1}
                    \draw (N\row\col) -- (N\nextrow\col);
                }
            }
        
            \node at (\cols/1.75, -\rows-0.5) {Grid M\"obius Strip: $\Mobius{m}{n}$};
            
        \end{tikzpicture}
        \end{center}
    \end{minipage}
    
    \end{center}

    \caption{The three types of grid-like graph we consider, each drawn with $m = 4$ and $n = 6$.}
    \label{gridlike example figure}
\end{figure}

    In Figure~\ref{gridlike example figure}, we depict elements from each of these families with $m=\rows$ and $n=\cols$. We note that while $\Band{m}{n}$ is isomorphic to $\Tube{n}{m}$, it will be natural for us to distinguish the two, as we will typically consider $m$ to be fixed and $n$ to be varying.

Just as for the rectangular grid-graphs, the behavior of $|\MIS(G)|$ is understandable for all grid-like graphs $G$; in a forthcoming paper, we give details on one possible enumeration scheme. In fact, the structure of grid-like graphs does not just make the enumeration problem tractable, but imparts significant structure to the graphs' set of MIS's. In this paper, we consider a variety of statistics of $\MIS(G)$ beyond its size. In particular, we are interested in the symmetries of grid-like graphs and their MIS's, and in the average size of MIS's in grid-like graphs. 

In a graph $G$, we say that two elements $I, I'$ of $\MIS(G)$ are \textit{isomorphic} if there is a graph automorphism $\xi: V(G) \rightarrow V(G)$ with $\xi(I) = I'$, and \textit{non-isomorphic} if no such $\xi$ exists. 
We let $\NIMIS(G)$ denote the set of isomorphism classes (equivalently, orbits under the action of the automorphism group) of MIS's in $G$. Non-isomorphic MIS's have been studied previously, for instance by Bisdorff and Marichal~\cite{Bisdorff} in the case where $G$ is a cycle. Ortiz and Villenueva~\cite{Ortiz} raised the question of counting non-isomorphic sets in grids and related graphs, which is naturally motivated:
the automorphisms of grid-like graphs are particularly appealing, as they naturally arise from the symmetries of their host surfaces. For flat rectangular grids (the case with fewest symmetries), we are able to exactly determine $|\NIMIS(\Grid{2}{n})|$ and to estimate $|\NIMIS(\Grid{m}{n})|$ in general.  Let $F(n)$ denote the $n$th Fibonacci number.

\begin{restatable}{thm}{twobyngrid}\label{2xn Grid Thm}
     For $n\geq 3$, we have
\[|\NIMIS(\Grid{2}{n})| = 
    \begin{cases}
    \frac{1}{2}\big(F(n)+F\big(\frac{n}{2}\big)\big) & \text{if $n$ is even, }\\
    \frac{1}{2}\big(F(n)+F\big(\frac{n+3}{2}\big)\big) & \text{if $n$ is odd. }\\
    \end{cases}\]
\end{restatable}

\begin{restatable}{thm}{gridnimisestimate}\label{Grid NIMIS estimate}
Fix an integer $m \geq 2$. Then
\[\underset{n \rightarrow \infty}{\lim} \frac{|\NIMIS(G_{m \times n})|}{|\MIS(G_{m \times n})|} = \frac{1}{4}.\]
\end{restatable}
We also exactly determine $|\NIMIS(\Tube{3}{n})|$.

\begin{restatable}{thm}{threebyntube}\label{3xn Tube Thm}
    For $n \ge 2$, we have \[|\NIMIS(\Tube{3}{n})| = 2^{n-3} + 2^{\lceil\frac{n-4}{2}\rceil}.\]
\end{restatable}
Our average size results are motivated by recent work on the average size of connected vertex sets in grids \cite{AndrewVince}, as well as previous work~\cite{Barbosa} on the possible sizes of MIS's in grid-like graphs. Let $A(G)$ denote the average size of an element of $\MIS(G)$. We find exact formulas for the average size of an MIS in $G_{2\times n}$, $\Band{2}{n}$, and $\Mobius{2}{n}$, whose asymptotic behavior can be cleanly described. 

\begin{restatable}{thm}{gridaveragesize}\label{Grid Average Size Thm}
We have
\[A(\Grid{2}{n}) = \left(\frac{\varphi}{\sqrt{5}}\right)n + \frac{2}{5} + o(1).\]\end{restatable}

\begin{restatable}{thm}{bandaveragesize}\label{Band Average Size Thm}
We have
\[A(\Band{2}{n}) = \left(\frac{\varphi}{\sqrt{5}}\right)n + o(1).\]
\end{restatable}

\begin{restatable}{thm}{mobiusaveragesize}\label{Mobius Average Size Thm}
We have
\[A(\Mobius{2}{n}) = \left(\frac{\varphi}{\sqrt{5}}\right)n + o(1).\]
\end{restatable}

The remainder of the paper is laid out as follows. In Section~\ref{defs}, we define some standard notation and collect some specific definitions which will be used throughout the paper. In Section~\ref{statistics}, we give proofs of our results. Section~\ref{parity} focuses on some preliminary results on the parity of $|\MIS(G)|$, Section~\ref{nimis subsection} focuses on determining $|\NIMIS(G)|$, including proofs of Theorems \ref{2xn Grid Thm}, \ref{Grid NIMIS estimate}, and~\ref{3xn Tube Thm}, and Section~\ref{avg size subsection} focuses on determining $A(G)$, including proofs of Theorems~\ref{Grid Average Size Thm}, \ref{Band Average Size Thm}, and~\ref{Mobius Average Size Thm}.

\subsection{Definitions and Notation}\label{defs}
Throughout the paper, we endeavor to use standard notation and terminology. We collect here some of the ubiquitous definitions and notational conventions for reference.

In a graph $G$, we denote the \textit{degree} of $u \in V(G)$ by $d(u)$.
For a positive integer $n$, we denote $\{1,2,\ldots,n\}$ by $[n]$. By $P_n$, we mean the graph with vertex set $[n]$ and edge set $\{\{i,j\} \in [n] : |i-j| = 1\}$, and by $C_n$, we mean $P_n$ with the edge $\{1,n\}$ added. Given two graphs $G$ and $H$, we denote the \textit{Cartesian product} of $G$ and $H$ by $G \Box H$.

A length-$n$ \textit{string} is a sequence of $n$ characters chosen from a pre-defined alphabet $\Sigma$, and we denote by $\Sigma^n$ the set of all possible length-$n$ strings over $\Sigma$. We call each character in a string a \textit{bit}. A length-$n$ binary string is a string from $\{0,1\}^n$. Given a string $s$, we use $s(i)$ to denote the $i$th bit of $s$.

We denote by $\mathbb{Z}_n$ the additive group of integers modulo $n$, and by $D_n$ the dihedral group of order $2n$, i.e., the group of symmetries of the regular $n$-gon.

$F(n)$ is the $n$th Fibonacci number, with the indexing convention $F(1) = F(2) = 1$, and $\varphi$ is the golden ratio.

\newcommand{\ver}{v}
\newcommand{\hor}{h}

\section{Proofs}\label{statistics}

\subsection{Symmetries}\label{parity}

Many of our arguments in  Sections~\ref{nimis subsection} and \ref{avg size subsection} will exploit the symmetries of grid-like graphs. We begin the study of these symmetries with a series of parity arguments, which introduce some ideas that will be useful in the upcoming subsections. 

\begin{defn} 
Let $G$ be an $m\times n$ grid-like graph. The \textit{$i$th slice} of $G$ is the subgraph of $G$ induced by the set of vertices with $x$-coordinate equal to $i$.
\end{defn}

For every grid-like graph $G$ except the M\"obius band, we have $G = H \Box P_n$ or $H \Box C_n$ for $H \in \{P_m,C_m\}$, and can equivalently say that the $i$th slice of $G$ is the subgraph of $G$ corresponding to the $i$th vertex of $P_n$ or $C_n$ in the Cartesian product. 
We can naturally think of grid-like graphs as arising from the concatenation of slices, and it can be convenient to consider MIS's in grid-like graphs by considering their intersection with individual slices. We begin with a simple observation about such intersections. Say that two slices $H_1,H_2$ of a grid-like graph $G$ are \textit{adjacent} if $G$ contains some edge between $V(H_1)$ and $V(H_2)$.

\begin{lemma}\label{lem:2emptycols}
    Let $H$ be a graph containing at least one edge and fix $n \geq 2$. Every MIS on $H \Box P_n$ or $H \Box C_n$ must contain at least one vertex from every pair of adjacent slices' connected components. 
\end{lemma}
\begin{proof}
    We consider $H \Box P_n$; the argument for $H \Box C_n$ is analogous. Suppose for a contradiction there exists $M \in \MIS(H \Box P_n)$ and adjacent slices  $H_1$ and $H_2$ such that $M$ is disjoint from $V(H_1) \cup V(H_2)$. If $H$ contains no slices distinct from $H_1,H_2$, then $M = \emptyset$, a contradiction to its maximality. So $H$ has at least $3$ slices, and thus at least one of $H_1,H_2$ (without loss of generality, $H_1$) is adjacent to a unique slice not in $\{H_1,H_2\}$. Let $H_0$ be this slice, and let $\iota_i : H \to H_i, i \in \{0,1,2\}$ be the embeddings defining each slice. Let $u_1$ and $u_2$ be two adjacent vertices in $H$ (which exist since $H$ contains at least one edge).

    Because $M$ is maximal, $\iota_1(u_1)$ 
    must be covered by some vertex. Then $\iota_0(u_1)$, the unique vertex in $H_0$ incident to $\iota_1(u_1)$, is the only other vertex which can cover $\iota_1(u_1)$, as no vertices in $H_1$ or $H_2$ are in $M$. So $\iota_0(u_1)$ must be in $M$. Similarly, $\iota_1(u_2)$ must be covered by $\iota_0(u_2)$, so $\iota_0(u_2)$ must be in $M$. However, $u_1,u_2 \in H$ are adjacent, so in the embedding $\iota_0(u_1), \iota_0(u_2) \in H_0$ are adjacent. This contradicts the assumption that $M$ is an independent set.
\end{proof}

By essentially the same argument, for $m \geq 2$, every $M \in \MIS(\Mobius{m}{n})$ intersects every pair of adjacent slices. (We simply need to verify that, although it cannot be described as a Cartesian product, $\Mobius{m}{n}$ still has the property that if $u_1,u_2 \in V(H_1)$ are adjacent, then their neighbors $\Tilde{u_1}, \Tilde{u_2}$ in an adjacent slice are also adjacent.) 

One application of Lemma~\ref{lem:2emptycols} is to determine the parity of $|\MIS(G)|$ for various grid-like graphs $G$. In order to express these parity arguments, we will also wish to formalize a notion of ``symmetric'' MIS's in grid-like graphs, which will also be useful later. 

\begin{defn}\label{symmetric MIS defn}
Suppose $G$ is a grid-like graph and $\xi: V(G) \rightarrow V(G)$ is an automorphism of $G$. We say that $M \in \MIS(G)$ is \textit{symmetric under} $\xi$ if $M = \xi(M)$.
\end{defn}

The description of MIS's as ``symmetric'' under automorphisms alludes to the geometric interpretation of automorphisms in this context. While the existence of symmetric MIS's is not tied to any particular depiction of their host graph, the grid-like graphs can be very naturally visualized as geometric objects (e.g., the grid-like graph $\Grid{m}{n}$ as a rectangular lattice) and their automorphisms as geometric symmetries of these objects. By considering which (if any) MIS's are symmetric under certain automorphisms, we can determine that almost all grid-like graphs have an even total number of MIS's.

Recall, $V_{m \times n} = \{(i,j) : 1 \leq j \leq m, 1 \leq i \leq n\}$ is the vertex set on which the $m \times n$ grid-like graphs are defined.  

\begin{defn}\label{def:automorphisms}
The \textit{vertical flip} $\ver: V_{m\times n} \to V_{m \times n}$ is $\ver(i,j) = (n+1 - i, j)$. The \textit{horizontal flip} $\hor: V_{m\times n} \to V_{m \times n}$ is $\hor(i,j) = (i, m+1 -j)$.
\end{defn}

Note that the vertical flip can be visualized as a reflection over a central vertical axis, and the horizontal flip can be visualized as a reflection over a central horizontal axis. Both are automorphisms for all grid-like graphs.

\begin{thm}\label{thm:grid-even} 
    $|\MIS(\Grid{m}{n})|$ is even for $m,n \ge 2$. 
\end{thm}

\begin{proof}
    
   We begin by directly proving the result for all even $n$. Assume $n = 2k$. 
    Consider some $M \in \MIS(\Grid{m}{n})$ which is not symmetric under the vertical flip $\ver$.
    Then $M, \ver(M)$ are distinct sets. It is straightforward to observe that $\ver(M) \in \MIS(\Grid{m}{n})$ and $\ver^2(M) = M$.  
    Thus, sets of the form $\{M, \ver(M)\}$ form a partition of these non-symmetric MIS's into sets of size $2$. Thus, the set of MIS's which are not symmetric under $\ver$ has even size. We now show that no MIS $M$ is symmetric under $\ver$. Assume (for a contradiction) that there exists such an $M$. Then for every $r \in [m]$, we have $(k,r) \in M$ if and only if $(k+1,r) \in M$. By the independence of $M$, we cannot have both $(k,r)$ and $(k+1,r)$ in $M$ for any $r$. Thus, no vertices in slices $k$ and $k+1$can be in $M$.
    By Lemma \ref{lem:2emptycols}, this contradicts the maximality of $M$. We conclude that $|\MIS(G_{m\times n})|$ is even when $n$ is even.

    For odd $n \geq 3$, we will proceed by (strong) induction on $n$. Fix $n = 2k+1 \geq 3$ and assume that for all $2 \leq \ell < n$, we have that $|\MIS(\Grid{m}{\ell})|$ is even. 
    Let $M \in \MIS(\Grid{m}{n})$. Note that the vertices in slice $k+1$ are fixed by $\ver$; consequently, it is possible that some MIS of $\Grid{m}{n}$ is symmetric under $\ver$.
    As in the even case, the set of MIS's which are not symmetric under $\ver$ is even. Now, we consider those MIS's which are symmetric under $\ver$. Suppose $M \in \MIS(\Grid{m}{n})$ has $M = \ver(M)$. We claim that the subset $M'$ of $M$ induced by the first $k+1$ slices of $G_{m\times n}$ is an MIS of $G_{m\times (k+1)}$. Certainly, $M'$ is independent. If it is not maximal, then there exists a vertex $u$ among the first $k+1$ slices which is uncovered by $M'$ but is covered by $M$. Since the vertices of $M\setminus M'$ are adjacent only to vertices in slices with index $k+1$ and higher, $u$ must be contained in slice $k+1$, and is covered by a vertex $w$ of $M$ in slice $k+2$. However, by the symmetry of $M$, if $w$ covers $u$, then a corresponding vertex $w'$ in slice $k$ covers $u$, a contradiction. Thus, $M'$ is a maximal independent set of $\Grid{m}{(k+1)}$. On the other hand, each MIS $M'$ of $\Grid{m}{(k+1)}$ can be used to form an MIS $M$ of $\Grid{m}{n}$ by setting $M = M' \cup \{(n+1-i,j): (i,j) \in M'\}$. (Geometrically, this corresponds to embedding $M'$ into the ``left half'' of the $m\times n$ rectangular grid and then reflecting $M'$ over the vertical axis to generate a vertically symmetric MIS of $\Grid{m}{n}$.) By strong induction, $|\MIS(\Grid{m}{(k+1)})|$ is even, and a bijection is created between $\MIS(\Grid{m}{(k+1)})$ and elements of $\MIS(\Grid{m}{n})$ which are symmetric under $\ver$. We conclude that $|\MIS(\Grid{m}{n})|$ is even.
\end{proof}

Similar ideas can be used to show that other grid-like graphs have an even number of MIS's. We begin with a result which implies that $|\MIS(\Band{m}{n})|$ and $|\MIS(\Tube{m}{n})|$ are always even. This result also applies to grid tori ($C_m \Box C_n$) as well as many other graphs not related to grids.

\begin{thm}\label{thm:HxCn-even} 
    Let $H$ be a graph with no isolated vertices. Then, $|\MIS(H \Box C_n)|$ is even for $n \ge 2$.
\end{thm}

\begin{proof}

     Consider some $M \in \MIS(H \Box C_n)$. Let $V(C_n) = [n]$ with the usual edges, let $H_i$ be the $i$th slice of $H \Box C_n$, and let $\iota_i$ be the embedding $H \to H_i$. Assume for the sake of contradiction that $M$ is symmetric under the vertical flip $\ver$. Then, assume there is some $u \in H$ such that $\iota_1(u) \in M$. But then by symmetry, we also have $\iota_n(u) \in M$, contradicting the assumption that $M$ is an independent set. Thus, no vertices in $H_{1}$ nor in $H_{n}$ belong to $M$, contradicting the maximality of $M$ by Lemma \ref{lem:2emptycols}. Thus $M$ is not symmetric under $\ver$, so $\{M, \ver(M)\}$ has cardinality $2$. Thus, $\MIS(H \Box C_n)$ can be partitioned into cardinality-2 sets, so it must have even cardinality.
\end{proof}

\begin{thm}\label{thm:mob-even}
    $|\MIS(\Mobius{m}{n})|$ is even for any $m,n \ge 2$.
\end{thm}

\begin{proof}

    Let $M \in \MIS(\Mobius{m}{n})$. If $M$ is not symmetric under the horizontal flip $\hor$, then $M$ $\{M, \hor(M) \}$ has cardinality $2$, so the set of non-symmetric MIS's is even. It thus suffices to show that the set of fixed points of $\MIS(\Mobius{m}{n})$ under $\hor$ is even. By Theorem \ref{thm:HxCn-even}, $|\MIS(\Band{m}{n})|$ is even, and the set of MIS's in $\MIS(\Band{m}{n})$ that are not symmetric under $\hor$ is even as always, so the set of fixed points under $\hor$ is even as well. We can then define a bijection of MIS's between the fixed points of $\MIS(\Mobius{m}{n})$ and the fixed points of $\MIS(\Band{m}{n})$ induced by the identity map on the set of vertices $V_{m \times n}$, as vertical symmetry ensures that an MIS is locally identical on either side of this correspondence.
\end{proof}

\subsection{Non-isomorphic MIS's}\label{nimis subsection}

The parity arguments in Section~\ref{parity} rely upon the correspondence of certain MIS's under some symmetry of a grid-like graph. When considering such symmetries, it is natural to ask how many MIS's are mutually distinct up to symmetry. Formally, in a graph $G$, we say that two elements $I, I'$ of $\MIS(G)$ are \textit{isomorphic} if there is a graph automorphism $\xi: V(G) \rightarrow V(G)$ with $\xi(I) = I'$, and \textit{non-isomorphic} if no such $\xi$ exists. To discuss non-isomorphic MIS's efficiently, we will want some standard terminology. Given a group $A$ and finite set $X$, we say that $A$ \textit{acts on} $X$ if there is a map $A\times X \rightarrow X$, denoted $(g,x) \mapsto g \cdot x$, satisfying two properties:
    \begin{enumerate}
        \item $e \cdot x = x$ for all $x \in X$, where $e$ is the identity element of $A$;
        \item $g_1 \cdot (g_2 \cdot x) = (g_1g_2) \cdot x$ for all $g_1,g_2 \in A$ and $x \in X$.
    \end{enumerate}
In particular, note that the automorphism group of a graph naturally acts on its set of MIS's, with $\xi \cdot M = \xi(M)$.

Given $A$ acting on $X$ and $x \in X$, we say that the \textit{orbit of $x$} is the set $\mathcal{O}(x) \coloneq \{g \cdot x | g \in A\}$, and the \textit{stabilizer of $x$} is the set $\Stab{x} \coloneq \{g \in A \ | \ g \cdot x = x\}$. Note that when $X = \MIS(G)$ and $A$ is the automorphism group of $G$, two MIS's  $M$ and $M'$ are isomorphic if and only if $\mathcal{O}(M) = \mathcal{O}(M')$. We define 
\[\NIMIS(G) := \{\mathcal{O}(M): M \in \MIS(G)\};\]
the goal of this section is to understand $|\NIMIS(G)|$ for certain grid-like graphs.

Using the phrasing of group actions will allow us to take advantage of pre-existing theory to concisely analyze $\NIMIS(G)$. For instance, the well-known Orbit-Stabilizer Theorem states that if $A$ acts on $X$ and $x \in X$, then
\[ |\mathcal{O}(x)| = \frac{|A|}{|\Stab{x}|};\]
in particular, we have that $|\mathcal{O}(x)|$ divides $|A|$, which will be repeatedly useful. 

We begin with consideration of non-isomorphic MIS's in grid graphs. Observe that if $m \neq n$, then $G_{m \times n}$ has four automorphisms, since under any automorphism $\xi$, the degree $2$ vertex $(1,1)$ must be mapped to one of the four degree $2$ vertices in $G_{m\times n}$, and the image of $(1,1)$ determines the images of every other vertex under $\xi$. (The automorphism group of the square grid graph $G_{n \times n}$ is $D_4$, which has order $8$; however, we will consider only non-square grid graphs here.) The four automorphisms of $\Grid{m}{n}$ are the identity map $id$, the vertical flip $\ver$, the horizontal flip $\hor$, and the double flip $\hor \circ \ver = \ver \circ \hor$, where $\ver$ and $\hor$ are as defined in Definition~\ref{def:automorphisms}.

We will first consider the $2 \times n$ grid graphs before generalizing to $m \times n$ grid graphs. For comparison, we note the total number of MIS's in $\Grid{2}{n}$. While $|\MIS(\Grid{m}{n})|$ is understandable for all $m$, the $m = 2$ case has a much simpler and more compact form than a general description of $|\MIS(\Grid{m}{n})|$. 
\begin{thm}[Euler~\cite{Euler2}]\label{2 by n count}
We have $|\MIS(\Grid{2}{n})| = 2F(n)$.
\end{thm}
Not only will Theorem~\ref{2 by n count} be quite relevant to our exploration of the statistics of $\MIS(\Grid{2}{n})$, but the value of $|\MIS(\Grid{2}{n})|$ suggests a connection to Fibonacci numbers which will be a recurring theme. 

The following function, which associates the MIS's in certain $2 \times n$ grid-like graphs with length-$n$ binary strings, will be a useful descriptive tool at several points.

 \begin{defn}\label{def:mappsi}
    Let $\psi:\MIS(\Grid{2}{n}) \rightarrow \{0,1\}^n$ be defined by \[\psi(M)(i) = |\{(i,1),(i,2)\} \cap M|\]
    for every $1 \le i \le n$ and every $M \in \MIS(\Grid{2}{n})$. 
 Let $\psi_c: \MIS(\Band{2}{n}) \cup \MIS(\Mobius{2}{n}) \rightarrow \{0,1\}^n$ be defined in the same way. For a length-$n$ string in the image of $\psi$, we say bits $i$ and $j$ are \emph{adjacent} if $|i-j|=1$. For a length-$n$ string in the image of $\psi_c$, which we will call a \emph{cyclic string}, bits $i$ and $j$ are adjacent if $|i-j| \in \{1,n-1\}$. 
 \end{defn}

Using $\psi$ will make it somewhat simpler to precisely determine $|\NIMIS(\Grid{2}{n})|$. We note the form which strings in the image of $\psi$ must take.

\begin{lemma}\label{map facts}
If $M \in \MIS(\Grid{2}{n})$, then the binary string $\psi(M)$ contains no two adjacent $0$'s, and its first and last bits are equal to $1$. Moreover, any length-$n$ binary string satisfying these conditions is the image of exactly two elements of $M, M' \in \MIS(\Grid{2}{n})$ under $\psi$, which satisfy $h(M) = M'$.
\end{lemma}

\begin{proof}
It follows immediately from Lemma~\ref{lem:2emptycols} that $\psi(M)$ contains no two adjacent $0$'s. Moreover, it is clear that for every $M \in \MIS(\Grid{2}{n})$, $M$ must intersect both the first and last slices of $\Grid{2}{n}$, so $\psi(M)$ must have first and last bits equal to $1$. 

To determine the pre-image of a length-$n$ binary string, we employ the following simple observation. Let $M \in \MIS(\Grid{2}{n})$ and suppose $j_1 < j_2 < \dots < j_k$ are the indices of slices of $\Grid{2}{n}$ which have non-empty intersection with $M$. Then either $(j_i, 1) \in M$ for all even $i$ and $(j_i, 2) \in M$ for all odd $i$, or vice versa. Indeed, suppose that this property does not hold. Then there exist $j_i, j_{i+1}$ such that either $(j_i, 1), (j_{i+1},1) \in M$ or $(j_i, 2), (j_{i+1},2) \in M$. Without loss of generality, say  $(j_i, 1), (j_{i+1},1) \in M$. Since $M$ is maximal, we have either $j_{i+1} = j_i + 1$ or $j_{i+1} = j_i + 2$ by Lemma~\ref{lem:2emptycols}. In the first case, $(j_i, 1), (j_{i}+1,1) \in M$ contradicts the independence of $M$. In the second, by the definition of the $j$ indices and independence of $M$, none of $(j_i,2), (j_{i} + 1, 1), (j_i + 1, 2), (j_i + 2, 2)$ are in $M$. This contradicts the maximality of $M$ (as $(j_i + 1, 2)$ could be added).

Using the above observation, note that each $n$-bit binary string can correspond to at most two MIS's, one including $(1,1)$ and one including $(1,2)$; moreover, if a binary string corresponds to two distinct MIS's $M$ and $M'$, then these will satisfy $M' = h(M)$. Thus, it suffices to show that each $n$-bit binary string without adjacent $0$'s and with starting and ending bits $1$ in fact corresponds to exactly two MIS's of $\Grid{2}{n}$. Given a binary string $b$ satisfying the stated properties, let $j_1, \dots, j_m$ be the indices of bits equal to $1$, let 
\[M_1 = \{(j_i, 1): \text{$i$ is even}\} \cup  \{(j_i, 2): \text{$i$ is odd}\}\]
and
\[M_2= \{(j_i, 2): \text{$i$ is odd}\} \cup  \{(j_i, 1): \text{$i$ is even}\}.\]
It is straightforward to observe that $M_1, M_2$ are independent sets of $\Grid{2}{n}$, so we must only show that each is maximal. For any $i$, each vertex in slice $j_i$ is either in $M_1$ or adjacent to an element of $M_1$. If slice $k$ has empty intersection with $M_1$, then, by the properties of $b$, we have $k \not \in \{1,n\}$, and slices $k-1,k+1$ intersect $M$. So either $(k-1,1), (k+1,2) \in M_1$ or  $(k-1,2), (k+1,1) \in M_1$; in either case, each vertex of slice $k$ is adjacent to a vertex of $M_1$. Thus, $M_1$ is maximal, and analogously, $M_2$ is maximal.
\end{proof}

\twobyngrid*

\begin{proof}

By Theorem~\ref{2 by n count}, we have $|\MIS(\Grid{2}{n})| = 2F(n)$. Thus, by Lemma~\ref{map facts}, 
\[|\psi(\MIS(\Grid{2}{n}))| = F(n).\]
Moreover, for any $M \in \MIS(\Grid{2}{n})$, we have $M \neq h(M)$, so $2 \leq |\mathcal{O}(M)|$, which implies $|\mathcal{O}(M)| \in \{2,4\}$ since the automorphism group of $\Grid{2}{n}$ has order $4$.  If elements of $\mathcal{O}(M)$ map to $k$ distinct strings under $\psi$, then by Lemma~\ref{map facts}, these correspond to $2k$ distinct MIS's which are pairwise isomorphic. Thus, we have that $|\psi(\mathcal{O}(M))| = \frac{|\mathcal{O}(M)|}{2} \in \{1,2\}$. So $|\psi(\MIS(\Grid{2}{n}))|$ counts orbits of size $2$ in $\NIMIS(\Grid{2}{n})$ once and orbits of size $4$ twice. To understand this double-counting, we characterize those strings corresponding to MIS's $M$ with $|\mathcal{O}(M)| = 2$. 

We shall say that an $n$-bit binary string $b_1b_2\dots b_n$  is \textit{vertically symmetric} if 
\[ b_1b_2\dots b_{n-1}b_n = b_n b_{n-1} \dots b_2 b_1.\]
We claim that the number of distinct orbits $\mathcal{O}(M)$ with $|\mathcal{O}(M)| = 2$ is precisely the number of vertically symmetric strings in the image of $\psi.$ Indeed, observe that if $\psi(M) = b_1b_2\dots b_n$, then $\psi(\ver(M)) = b_n \dots b_2b_1$. Thus, if $\psi(M)$ is vertically symmetric, then we have $\psi(\ver(M)) = \psi(M)$, so $\ver(M) \in \{M, \hor(M)\}$, which implies $|\mathcal{O}(M)| = 2$. On the other hand, if $\psi(M)$ is not vertically symmetric, then $\psi(\ver(M)) \neq \psi(M)$, so $v(M) \not \in \{M,\hor(M)\}$, which implies $|\mathcal{O}(M)| =  4$.

Let $VS$ be the set of vertically symmetric strings in the image of $\psi$. We thus have 
\[|\psi(\Grid{2}{n})| + |VS| = 2|\NIMIS(\Grid{2}{n})|.\]
Using the known value of $|\psi(\Grid{2}{n})|$, this yields
\[\frac{F(n) + |VS|}{2} = |\NIMIS(\Grid{2}{n})|.\]
To finish, we must determine $|VS|$. We count the number of vertically symmetric binary strings of length $n$ which satisfy the conditions of Lemma~\ref{map facts}. Suppose $b_1b_2 \dots b_n$ is a vertically symmetric string in the image of $\psi$. There are two cases.

\textbf{Case 1: $n=2k$ for some $k \ge 2$.} Here, the line of symmetry falls between the $k$th and $(k + 1)$st bits. That is, $b_k$ must equal $b_{k+1}$. By Lemma~\ref{lem:2emptycols}, $b_{k} = b_{k + 1} = 1$. Thus, $b_1b_2\dots b_k$ is a binary string of length $k$ which starts and ends in a $1$ and has no adjacent $0$'s. The number of such strings is $|\psi(\Grid{2}{k})| = F(k)$ by Lemma~\ref{map facts}. Observe also that an element of $\psi(\Grid{2}{k})$ uniquely determines a vertically symmetric element of $\psi(\Grid{2}{n})$ (by reflection over the vertical line of symmetry). So we have $|VS| = F(k) = F(\frac{n}{2})$. 

\textbf{Case 2: $n=2k+1$ for some $k \ge 1$.} Here, the axis of symmetry falls on the  $(k+1)$st bit. If $b_{k+1} = 1$, then $b_1b_2 \dots b_{k+1}$ is an element of $\psi(\Grid{2}{(k+1)})$. By reflection, each such element uniquely determines a vertically symmetric element of $\psi(\Grid{2}{n})$ with $b_{k+1} = 1$, so there are $F(k+1)$ vertically symmetric elements of $\psi(\Grid{2}{n})$ with $b_{k+1} = 1$.
    If $b_{k+1} = 0$, then by Lemma~\ref{lem:2emptycols}, we must have $b_{k} = 1$. Thus, $b_1 b_2 \dots b_{k}$ is an element of $\psi(\Grid{2}{k})$. By reflection, each such element uniquely determines a vertically symmetric element of $\psi(\Grid{2}{n})$ with $b_{k+1} = 0$, so there are $F(k)$ vertically symmetric elements of $\psi(\Grid{2}{n})$ with $b_{k+1} = 0$. In total, $|VS| = F(k+1) + F(k) = F(k+2) = F\big(\frac{n+3}{2}\big)$.

We conclude that when $n$ is even,
\[|\NIMIS(\Grid{2}{n})| = \frac{1}{2}\left(F(n)+F\left(\frac{n}{2}\right)\right),\]
and when $n$ is odd,
\[|\NIMIS(\Grid{2}{n})| = \frac{1}{2}\left(F(n) + F\left(\frac{n+3}{2}\right)\right).\]
\end{proof}

The use of $\psi$ somewhat simplifies the argument to determine $|\NIMIS(\Grid{2}{n})|$ because $\psi$ ``compresses'' the relevant information about each MIS of $\Grid{2}{n}$ into a form which is somewhat easier to manage than the full MIS. By working with $\psi(M)$ instead of $M$, we can bypass describing exactly when $v(M) = M$ and when $(v \circ h)(M) = M$. A direct proof of the enumeration, without $\psi$, is possible but somewhat more cumbersome; to explicitly describe the construction of MIS's with orbits of size $2$, it becomes necessary to pinpoint when it is possible for each isomorphism to fix $M$.

For larger $m$, an exact determination of $|\NIMIS(\Grid{m}{n})|$ seems more difficult. As the possible slices of an MIS become more numerous and varied, it becomes more complicated to describe when a ``partial'' independent set can be reflected over an axis of symmetry to actually create an MIS in the full grid. However, the fact that each $M \in \MIS(\Grid{m}{n})$ with $|\mathcal{O}(M)| < 4$ is in some way determined by only ``half'' of the grid will allow us to asymptotically estimate $|\NIMIS(\Grid{m}{n})|$. Note that, combining Theorems~\ref{2 by n count} and \ref{2xn Grid Thm}, we have that $|\NIMIS(\Grid{2}{n})|$ can be approximated by $\frac{F(n)}{2} = \frac{|\MIS(\Grid{2}{n})|}{4}$, as $F(n)$ dominates $F\left(\frac{n}{2}\right)$ and $F\left(\frac{n+3}{2}\right)$. This growth behavior corresponds to the observation that, as $n$ grows, the proportion of binary strings which are vertically symmetric (and thus, of MIS's of $\Grid{2}{n}$ which are symmetric under some automorphism of $\Grid{2}{n}$) goes to zero. The following theorem generalizes this observation.

\gridnimisestimate*

\begin{proof}

We fix  $m,n\in \mathbb{N}$. Since we will consider a limit as $n$ approaches infinity, we may assume $m < n$; in particular, the automorphism group of $\Grid{m}{n}$ has order $4$. We will also assume $\lceil \frac{n}{2} \rceil - 4 > 0$. Note that when $m = 2$, the result is immediately implied by Theorem~\ref{2xn Grid Thm}, so we may assume $m \geq 3$. Let 
\[\mathcal{S}_i := \{\mathcal{O}(M) : M \in \MIS(\Grid{m}{n}), |\mathcal{O}(M)| = i.\}\] 
Then, since $\NIMIS(\Grid{m}{n})=\mathcal{S}_4 \cup \mathcal{S}_2 \cup \mathcal{S}_1$, we have

\[4 \cdot \frac{|\NIMIS(\Grid{m}{n})|}{|\MIS(\Grid{m}{n})|} = \frac{4|\mathcal{S}_4| + 4|\mathcal{S}_{2}| + 4|\mathcal{S}_1|}{|\MIS(\Grid{m}{n})|}.\]
Note that 
$4|\mathcal{S}_4| + 2|\mathcal{S}_{2}| + |\mathcal{S}_1| = |\MIS(\Grid{m}{n})|,$ 
so we have
\[1 = \frac{4|\mathcal{S}_4| + 2|\mathcal{S}_{2}| + |\mathcal{S}_1|}{|\MIS(\Grid{m}{n})|} \leq 4 \cdot \frac{|\NIMIS(\Grid{m}{n})|}{|\MIS(\Grid{m}{n})|}\]
and so $\frac{|\NIMIS(\Grid{m}{n})|}{|\MIS(\Grid{m}{n})|}$ is bounded below by $\frac14$ for any fixed $n$. Also, we have
\[4 \cdot \frac{|\NIMIS(\Grid{m}{n})|}{|\MIS(\Grid{m}{n})|} \leq   \frac{4|\mathcal{S}_4| + 2|\mathcal{S}_{2}| + |\mathcal{S}_1|}{|\MIS(\Grid{m}{n})|} +\frac{ 3\big(|\mathcal{S}_{2}| + |\mathcal{S}_1|\big)}{|\MIS(\Grid{m}{n})|} = 1 + \frac{ 3\big(|\mathcal{S}_{2}| + |\mathcal{S}_1|\big)}{|\MIS(\Grid{m}{n})|}. \]
and so $\frac{|\NIMIS(\Grid{m}{n})|}{|\MIS(\Grid{m}{n})|}$ is bounded above by $\frac14 + \frac{3\left(|S_2| + |S_1|\right)}{4|\MIS(\Grid{m}{n})|}$ for any fixed $n$. Thus, to prove the desired limit, it will suffice to show that \[\underset{n\rightarrow \infty}{\lim} \frac{ 3\big(|\mathcal{S}_{2}| + |\mathcal{S}_1|\big)}{4|\MIS(\Grid{m}{n})|} = 0.\]

To do so, we give a (rough) estimate of $|\mathcal{S}_2| + |\mathcal{S}_1|$. 
Clearly, if $\mathcal{O}(M) \in \mathcal{S}_2 \cup \mathcal{S}_1$, then $M$ is equal to (at least) one of $h(M), v(M), (v\circ h)(M)$. 
If $M = v(M)$ or $M = (v \circ h)(M)$, note that $M$ is determined by its intersection with the first $n' := \lceil \frac{n}{2} \rceil$ slices of $\Grid{m}{n}$. 
As in the proof of Theorem~\ref{thm:grid-even}, if $M = v(M)$, then the subset of $M$ induced by the first $n'$ slices of $\Grid{m}{n}$ is an element of $\MIS(\Grid{m}{n'})$, so there are at most $|\MIS(\Grid{m}{n'})|$ elements $M \in \MIS(\Grid{m}{n})$ with $M = v(M)$. 
If $M = (v \circ h)(M)$, then it is no longer necessarily true that the subset $M'$ of $M$ induced by the first $n'$ slices of $\Grid{m}{n}$ is an element of $\MIS(\Grid{m}{n'})$. 
However, we can map $M'$ to an element $M_+^{'} \in  \MIS(\Grid{m}{(n'+1)})$ as follows. Consider the independent set $I$ of $\Grid{m}{(n' + 1)}$ induced by the first $n' + 1$ slices of $M$. If $I$ is maximal, set $M_+' : = I$. If not, then (since $M$ is maximal), the only uncovered vertices in $I$ are in its $(n' + 1)$st slice, and we may obtain a maximal independent set $M_+'$ by adding some set of vertices in the $(n' + 1)$st slice to $I$. Thus, there exists (at least) one element of $\MIS\big(\Grid{m}{(n' + 1)}\big)$ whose first $n'$ slices induce $M'$. It follows that there are at most $\left|\MIS\big(\Grid{m}{(n' + 1)}\big)\right|$ elements $M \in \MIS(\Grid{m}{n})$ with $M = (v \circ h)(M)$. Finally, let $m' := \lceil \frac{m}{2} \rceil$: if $M = h(M)$, then $M$ is determined by its intersection with the vertices 
\[V' := \{(i,j): 1 \leq i \leq n, 1 \leq j \leq m' \}.\] 
Similarly to the proof of Theorem~\ref{thm:grid-even}, the subset $M \cap V'$ of $M$ is a maximal independent set of $\Grid{m'}{n}$. Thus, there are at most $|\MIS(\Grid{m'}{n})|$ elements $M \in \MIS(\Grid{m}{n})$ with $M =  h(M)$. Combining these bounds, we can write

\[\frac{ 3\big(|\mathcal{S}_{2}| + |\mathcal{S}_1|\big)}{4|\MIS(\Grid{m}{n})|} \leq \frac{ 3\big(|\MIS(\Grid{m}{n'})| + \left|\MIS(\Grid{m}{\left(n'+1\right)})\right| + |\MIS(\Grid{m'}{n})|\big)}{4|\MIS(\Grid{m}{n})|}. 
\]
We now bound each term in the numerator. We note that similar but less explicit bounds on these terms can be obtained from existing theory regarding the growth-behavior of $|\MIS(\Grid{m}{n})|$; however, in the interest of self-containment, we present direct arguments when possible.

Firstly, since $|\MIS(\Grid{m}{n'})| \leq \left|\MIS\big(\Grid{m}{(n' + 1)}\big)\right|$, we bound both of the first terms by considering $\left|\MIS\big(\Grid{m}{(n' + 1)}\big)\right|$. We can create an element $M \in \MIS(\Grid{m}{n})$ as follows. Choose an element $M_1 \in \MIS\big(\Grid{m}{(n' + 1)}\big)$ and an element $M_2 \in \MIS\big(\Grid{m}{(n' - 3)}\big)$. Let $M^-$ be the independent set of $\Grid{m}{n}$ whose first $n' + 1$ slices induce $M_1$, whose last $n' - 3$ slices induce $M_2$, and which contains no vertex which is not either in the first $n' + 1$ or the last $n' - 3$ slices. Note that since $$(n'+1) + (n'-3) = \left(\left\lceil \frac{n}{2} \right\rceil + 1\right) +  \left(\left\lceil \frac{n}{2} \right\rceil -3\right) \leq n + 1 + 1 - 3 = n-1,$$ $M^-$ contains at least one empty slice between the slices inducing $M_1$ and the slices inducing $M_2$. Thus, $M^-$ is in fact independent. We set $M := M^-$ if $M^-$ is maximal. If $M^-$ is not maximal, note than $M^-$ is contained in some $M \in \MIS(\Grid{m}{n})$ which (by the maximality of $M_1,M_2$) induces $M_1$ on its first  $n' + 1$ and $M_2$ on its last $n' - 3$ slices. In particular, we have
\[\left|\MIS\big(\Grid{m}{(n' + 1)}\big)\right| \cdot \left|\MIS\big(\Grid{m}{(n' - 3)}\big)\right| \leq |\MIS(\Grid{m}{n})|,\]
or equivalently,
\[\left|\MIS\big(\Grid{m}{(n' + 1)}\big)\right|  \leq \frac{|\MIS(\Grid{m}{n})|}{\left|\MIS\big(\Grid{m}{(n' - 3)}\big)\right|}.\]
Analogously, if $m' + 1 < m$, we can find an injection from $\MIS(\Grid{m'}{n}) \times \MIS(\Grid{1}{n})$ into $\MIS(\Grid{m}{n})$ to bound 
\[|\MIS(\Grid{m'}{n})| \leq \frac{|\MIS(\Grid{m}{n})|}{|\MIS(\Grid{1}{n})|}.\] 
If $m'+1=\lceil \frac{m}{2} \rceil +1 \geq m$, then $m \leq 3$. Recall that we may assume $m \geq 3$, so $m = 3$ is the only outstanding case, which we treat in a moment. For $m \geq 4$, set 
\[f_m(n) := \max\left\{ \frac{1}{\left|\MIS\big(\Grid{m}{(n'-3)}\big)\right|}, \frac{1}{|\MIS(\Grid{1}{n})|} \right\}.\]
Note that by our previous bounds, 
\[\frac{ 3\big(|\mathcal{S}_{2}| + |\mathcal{S}_1|\big)}{4|\MIS(\Grid{m}{n})|} \leq \frac{ 3\big(|\MIS(\Grid{m}{n'})| + \left|\MIS(\Grid{m}{\left(n'+1\right)})\right| + |\MIS(\Grid{m'}{n})|\big)}{4|\MIS(\Grid{m}{n})|} \leq \frac{9f_m(n)}{4}.\]
In particular, it is clear that $\underset{n\rightarrow \infty}{\lim} f_m(n) = 0$, so $\underset{n\rightarrow \infty}{\lim} \frac{ 3\big(|\mathcal{S}_{2}| + |\mathcal{S}_1|\big)}{4|\MIS(\Grid{m}{n})|} = 0$, as desired. When $m = 3$, the argument is the same, but we set 
\[f_m(n) := \max\left\{ \frac{1}{\left|\MIS\big(\Grid{m}{(n' -3 )}\big)\right|}, \frac{|\MIS(\Grid{m'}{n})|}{|\MIS(\Grid{m}{n})|} \right\}.\]
It remains true that 
\[\frac{ 3\big(|\mathcal{S}_{2}| + |\mathcal{S}_1|\big)}{4|\MIS(\Grid{m}{n})|} \leq \frac{9f_m(n)}{4},\]
and to show that $\underset{n\rightarrow \infty}{\lim} f_m(n) = 0$, it suffices to show that $$\lim_{n\rightarrow \infty} \frac{|\MIS(\Grid{m'}{n})|}{|\MIS(\Grid{m}{n})|}=\lim_{n\rightarrow \infty} \frac{|\MIS(\Grid{2}{n})|}{|\MIS(\Grid{3}{n})|} = 0.$$ This follows from the known asymptotic growth behavior of $|\MIS(\Grid{2}{n})|$ and $|\MIS(\Grid{3}{n})|$. By Theorem~\ref{2 by n count}, $|\MIS(\Grid{2}{n})| = \Theta(\varphi^n)$, while $|\MIS(\Grid{3}{n})|$ grows with a strictly larger exponential base (approximately 2.037); see \cite{Euler2, BigOh, Ortiz} for enumerative results from which this base may be inferred. 
\end{proof}

\newcommand{\rot}{r}
\newcommand{\refl}{h}
\newcommand{\vflip}{v}

We now turn to enumeration of $\NIMIS(G)$ for other grid-like graphs $G$, beginning with $\Tube{3}{n}$. We start by describing the general form of the automorphisms of $\Tube{m}{n}$. Consider a vertex $u$ contained in slice $1$ of $\Tube{m}{n}$. Then $d(u) = 3$. Observe that any automorphism $\xi$ of $\Tube{m}{n}$ must map $u$ to a vertex in either slice $1$ or slice $n$ of $\Tube{m}{n}$ in order to have $d(\xi(u)) = 3$. 
In order to preserve adjacencies under $\xi$, the mapping of slice $1$ into slice $i$ (for $i \in \{1,n\})$ must correspond to an element of the dihedral group $D_m$. 
Finally, note that the image of slice $1$ under $\xi$ uniquely determines $\xi$. Precisely, the automorphism group of $\Tube{m}{n}$ is (isomorphic to) $\mathbb{Z}_2 \times D_m$. In particular, the automorphism group of $\Tube{3}{n}$ is isomorphic to $\mathbb{Z}_2 \times D_3$, so
has order $12$ and is generated by the following elements, labeled according to the convention of Definition~\ref{def:automorphisms}. (We remark that in the context of the dihedral group, the horizontal flip is often called the reflection.)

\begin{itemize}
    \item The rotation $\rot$ determined by $\rot(1,1)=(1,2)$, $\rot(1,2)=(1,3)$, and $\rot(1,3)=(1,1)$;
    \item The horizontal flip $\refl$ determined by $\refl(1,1)=(1,3)$, $\refl(1,2)=(1,2)$, and $\refl(1,3)=(1,1)$;
    \item The vertical flip $\vflip$ determined by $\vflip(1,1)=(n,1)$, $\vflip(1,2)=(n,2)$, and $\vflip(1,3)=(n,3)$.
\end{itemize}
With these automorphisms described and labeled, we are prepared to count $|\NIMIS(\Tube{3}{n})|$.

\threebyntube*

\begin{proof}
    We begin by determining the possible sizes of an orbit $\mathcal{O}(M)$. First, note that for any $M \in\MIS(\Tube{3}{n})$ and any $i \in  [n]$, $M$ intersects the $i$th slice of $\Tube{3}{n}$ in exactly one vertex. In particular, $M$ intersects the $1$st slice in exactly one vertex, as do $\rot(M)$ and $\rot^2(M)$. Specifically, it is clear that $M, \rot(M),$ and $\rot^2(M)$ in fact intersect the $1$st slice in pairwise distinct vertices and so are pairwise distinct MIS's. Let $S_M := \{M, \rot(M), \rot^2(M)\}$ and
    
    \[\mathcal{P} := \{S_M : M \in \MIS(\Tube{3}{n})\}\]
    
    It is clear that $\mathcal{P}$ is a partition of $\MIS(\Tube{3}{n})$ into sets of size $3$. Moreover, the elements of a part are pairwise isomorphic. This immediately implies that every $M$ has $|\mathcal{O}(M)| \geq 3$; by considering classes in $\mathcal{P}$, we shall be able to improve this bound.

    Firstly, we observe that for any $S \in \mathcal{P}$ and for any $M \in S$, we have $\refl(M) \not\in S$. Indeed, let $M_1$ be the unique element of $S$ containing $(1,2)$. Observe that exactly one of $(2,1), (2,3) \in M_1$. Whichever of these two vertices $M_1$ contains, $\refl(M_1)$ contains the other, so $h(M_1) \neq M_1$. Since $(1,2) \in h(M_1)$, this implies that $\refl(M_1)$ and $M_1$ are in different parts of $\mathcal{P}$. Next, suppose $\refl(\rot^k(M_1)) \in S$ for some $k \in \{1,2\}$. Then, since $\refl \rot^k = \rot^{3-k}\refl$, we have $\rot^{3-k}(\refl(M_1)) \in S$. Because $S$ is closed under the action of $\rot$, this implies that $\rot^3 (\refl(M_1)) = \refl(M_1) \in S$, a contradiction.  In fact, this argument implies that $\refl$ maps the elements of $S$ to the elements of precisely one set in $\mathcal{P}$, say $S'$. Since $\refl$ has order $2$, it is also clear that $\refl$ maps the elements of $S'$ to the elements of $S$; thus, $\refl$ further partitions the sets in $\mathcal{P}$ into pairs $(S,S')$ such that the six distinct elements of $S \cup S'$ are pairwise isomorphic. This implies that in fact, every $M$ has $|\mathcal{O}(M)| \geq 6$, so every orbit has size either $6$ or $12$. Let $\mathcal{S}_6$ be the set of parts in $\mathcal{P}$ which contain MIS's with orbit size $6$. Orbits of size $12$ correspond to sets of $4$ parts from $\mathcal{P} \setminus \mathcal{S}_6$, while orbits of size $6$ correspond to sets of $2$ parts from $\mathcal{S}_6$. We thus have   
    \[\frac{|\mathcal{P}| + |\mathcal{S}_6|}{4} = |\NIMIS(\Tube{3}{n})|.\] 
    So, our goal is to determine $|\mathcal{P}|$ and $|\mathcal{S}_6|$. As in the proof of Theorem~\ref{2xn Grid Thm}, we begin by associating MIS's with particular strings which will clarify the action of the automorphisms on these MIS's.

    Let $I(M,i)$ be the $y$-coordinate of the unique vertex in the intersection between an MIS $M$ and the $i$th slice of $\Tube{3}{n}$. We define 
    \[\psi: \MIS(\Tube{3}{n}) \rightarrow \{+,-\}^{n-1}\]
    by setting $\psi(M)$ equal to the $(n-1)$-bit string with
    
    \[\psi(M)(i) := \begin{cases}
    + & \text{if $I(M,i+1)  \equiv I(M,i) + 1 \Mod 3$,} 
    \\  - & \text{if $I(M,i+1)  \equiv I(M,i) - 1 \Mod 3$}. 
    \end{cases}\]    

    We begin with some observations about the images of MIS's under $\psi$. Firstly, note that $\psi$ establishes a bijection between the parts of $\mathcal{P}$ and $\{+,-\}^{n-1}$. Indeed, $\psi(M) = \psi(M')$ if and only if $M$ and $M'$ are in the same class of $\mathcal{P}$, and given a bit-string $b \in \{+,-\}^{n-1}$, it is clear how to construct an MIS of $\Tube{3}{n}$ with image $b$ under $\psi$. Moreover, $\psi(M)$ is closely related to $\psi(\vflip(M))$ and $\psi(h(M))$: $\psi(h(M))$ is obtained from $\psi(M)$ by negating each sign of $\psi(M)$, while $\psi(\vflip(M))$ is obtained from $\psi(M)$ by reversing $\psi(M)$ and then negating each sign.

    We now use the correspondence between parts of $\mathcal{P}$ and these strings to count orbit sizes. Suppose $M \in S \in \mathcal{P}$, and $S'$ is the part of $\mathcal{P}$ containing $h(M)$. Note that $|\mathcal{O}(M)| = 12$ if and only if $\vflip(M) \not \in S \cup S'$. If $\vflip(M) \in S \cup S'$ (i.e., if $|\mathcal{O}(M)| = 6$), then either $\psi(\vflip(M)) = \psi(M)$  or $\psi(\vflip(M)) = \psi(\refl(M))$. We count the number of strings in $\{+,-\}^{n-1}$ which can satisfy either equality, breaking into cases according to parity of $n$. 

    \textbf{Case 1: $n = 2k+1$ for some $k\geq 2$.} In this case, $\psi(M)$ is of even length. If $\psi(M) = \psi(v(M))$, then $\psi(M)$ is self-equal when reversed and negated. Such a string is uniquely determined by its first $k$ bits (its $(n - i)$th bit must be equal to the negation of its $i$th bit for $1 \leq i \leq k$).
    These first $k$ bits can be fixed in $2^{k}$ ways. If $\psi(\vflip(M)) = \psi(h(M))$, then the negation of $\psi(M)$ is equal to the reversal and negation of $\psi(M)$. Again, a string with this property is uniquely determined by its first $k$ bits (its $(n - i)$th bit must be equal to its $i$th bit for $1 \leq i \leq k$).
    These first $k$ bits can be fixed in $2^{k}$ ways. Thus, we have $|\mathcal{P}|  =  |\{+,-\}^{n-1}| = 2^{n-1}$ and $|\mathcal{S}_6| = 2 \cdot 2^k = 2 \cdot 2^\frac{n-1}{2} = 2^{\frac{n+1}{2}}$, so
    \[|\NIMIS(\Tube{3}{n})| = \frac{1}{4}\left(2^{n-1} + 2^{\frac{n+1}{2}}\right) = 2^{n-3} + 2^{\frac{n-3}{2}}.\]

    \textbf{Case 2: $n=2k$ for some $k \geq 2$.} In this case, $\psi(M)$ is of odd length. Thus, when $\psi(M)$ is reversed, the $k$th bit is fixed. If $\psi(v(M)) = \psi(M)$, then $\psi(M)$ is self-equal under reversal and negation; this is impossible, as the $k$th bit cannot be both $+$ and $-$. If $\psi(\vflip(M)) = \psi(\refl(M))$, then the negation of $\psi(M)$ is equal to the reversal and negation of $\psi(M)$. A string with this property is unique determined by its first $k$ bits (its $(n-i)$th bit must equal its $i$th bit, for $1 \leq i \leq k - 1$). Thus, we have $|\mathcal{P}|  =  |\{+,-\}^{n-1}| = 2^{n-1}$ and $|\mathcal{S}_6| = 2^k = 2^{\frac{n}{2}}$, so \[|\NIMIS(\Tube{3}{n})| = \frac{1}{4}\left(2^{n-1} + 2^{\frac{n}{2}}\right) = 2^{n-3} + 2^{\frac{n-4}{2}}.\] Lastly, the cases of $n\in\{2,3\}$ can quickly be verified. The equality is thus proven.  
    \end{proof}

Finally, we consider non-isomorphic MIS's in the $2 \times n$ fat cylinder.
The automorphism group of $\Band{2}{n}$ is generated by $3$ elements:
\begin{itemize}
    \item The rotation $r$ determined by $r(i,1) = (i+1,1)$ for all $i \in [n-1]$ and $r(n,1) = (1,1)$; 
    \item The vertical flip $v$ determined by $v(i,1) = (n-i+1,1)$ for all $i \in [n]$;
    \item The horizontal flip $h$ determined by $h(i,1) = (i,2)$ and $h(i,2) = (i,1)$ for all $i \in [n]$. 
\end{itemize}
Recall the map $\psi_c$ described in Definition~\ref{def:mappsi}. Similar to the proof of Lemma~\ref{map facts}, we observe the following facts. The restriction that the total number of 1's in each string be even corresponds to the fact, previously observed by Barbosa and Cappelle (Lemma 4 in~\cite{Barbosa}), that every MIS of $\Band{2}{n}$ has even size. Since this parity condition will be relevant both here and in Section~\ref{avg size subsection}, we briefly indicate a proof in the interest of self-containment.

\begin{obs}\label{cyclic map facts}
If $M \in \MIS(\Band{2}{n})$, then the cyclic string $\psi_c(M)$ contains no two adjacent $0$'s (recall that the first and last bits are adjacent) and contains an even total number of $1$'s.
Moreover, any length-$n$ binary string satisfying these conditions is the image of exactly two elements $M,M' \in \MIS(\Band{2}{n})$ under $\psi_c$, which satisfy $h(M) = M'$.
\end{obs}

\begin{proof}[Proof sketch]
For all observations except that each cyclic string $\psi_c(M)$ contains an even number of $1$'s, an analogous argument to the proof of Lemma~\ref{map facts} suffices. We now discuss the parity condition. Suppose $M \in \MIS(\Band{m}{n})$ and say $j_1 < j_2 < \dots < j_k$ are the indices where $\psi_c(M)$ is equal to $1$. As in the proof of Lemma~\ref{map facts}, either $(j_i,1) \in M$ for all odd $i$ and $(j_i,2) \in M$ for all even $i$, or vice versa. Thus, if $k$ is odd, then without loss of generality, both $(j_1,1),(j_k,1) \in M$. This contradicts the independence of $M$ if slices $j_1,j_k$ are adjacent and contradicts the maximality of $M$ otherwise (since then the vertex $(j_k+1,2)$ could be added to $M$).
\end{proof}

We are able to describe $|\NIMIS(\Band{2}{n})|$, albeit with a less explicit formula. Recall that for a group $G$ acting on a set $X$, the set of orbits of this action is denoted $X/G$.
Additionally, recall that a $k$-composition of $n$ is an ordered sequence of $k$ strictly positive integers that sum to $n$.

\begin{thm}
    Let $C_k(n)$ be the set of $k$-compositions of $n$. Encode each composition in $C_k(n)$ as a cyclic string of $k$ integers adding to $n$, where cyclic rotations and reflections of the string are considered symmetries.
    Let $D_k$ denote the group of symmetries of a $k$-gon. Then, for all $n\geq3$, we have \[|\NIMIS(\Band{2}{n})| = \begin{cases}
        1 + \sum\limits^{\lfloor\frac{n}{4}\rfloor}_{k = 1}|C_{2k}(n-2k)/D_{2k}| &  \text{ if $n$ is even,}\\
        1 + \sum\limits^{\lfloor\frac{n-2}{4}\rfloor}_{k = 1}|C_{2k+1}(n-2k-1)/D_{2k+1}| & \text{ if $n$ is odd.}
    \end{cases}\]
\end{thm}
\begin{proof} 
For a set of strings $L$, we define the \textit{string reflection} 
    $V_L: L \rightarrow L$ by \[V_L(b_1b_2\dots b_{n-1}b_n) = b_nb_{n-1}\dots b_2b_1\] 
    and the \textit{string rotation} 
    $R_L: L \rightarrow L$ by \[R_L(b_1b_2\dots b_{n-1}b_n) = b_nb_1b_2 \dots b_{n-1}.\] 
    Let $A := \psi_c(\MIS(\Band{2}{n}))$. Using Observation~\ref{cyclic map facts} and the observation that $\psi_c(v(M)) = V_A(\psi_c(M))$ and $\psi_c(r(M)) = R_A(\psi_c(M))$,
    counting $|\NIMIS(\Band{2}{n})|$ translates to counting the number of distinguishable elements of $A$ with respect to string rotation and reflection.

    First, assume $n$ is even. Then from Observation~\ref{cyclic map facts}, each string in $A$ contains an even number of 1's and of 0's. Let $k>0$ and let $X_{2k}$ be the set of cyclic length-$n$ strings of $A$ with $(2k)$-many 0's. We show that $|X_{2k}/D_{n}| = |C_{2k}(n-2k)/D_{2k}|$.

    Because strings in $X_{2k}$ are cyclic, every orbit in $X_{2k}/D_{n}$ contains at least one representative with a $0$ at index $1$. Define $X'_{2k}:=\{x \in X_{2k} : x(1)=0\}$. It follows that there is at least one element in $X'_{2k}$ from each element of $X_{2k}/D_n$.

    There is a natural bijection $f: X'_{2k} \to C_{2k}(n-2k)$ defined as follows: Let $x \in X'_{2k}$. Then, $x(1) = 0$ and $x$ is composed of $(2k)$-many runs of $1$'s, each separated by exactly one $0$. So, set $f(x)$ to be the $(2k)$-composition of $n-2k$ realized by setting the $i$th integer of the composition to be length of the $i$th run of $1$'s in $x$. Clearly, $f$ is injective, and by Observation~\ref{cyclic map facts}, this map is surjective.

    Now, suppose $x, x' \in X'_{2k}$ are in the same orbit of $X_{2k}$ under $D_n$. Suppose $f(x) = (x_1, x_2, \dots, x_{2k})$. Since $x, x'$ are in the same orbit, $x'$ may be obtained from $x$ by a series of rotations and possibly a reflection. If no reflection is needed, then for some $i$,
    $$f(x') = (x_{2k-i}, x_{2k-i+1}, \dots, x_{2k}, x_1, \dots, x_{2k-i-1}).$$
    In this case, $f(x)$ is in the same orbit as $f(x')$ under $D_{2k}$ (rotate $f(x)$ $(i+1)$-many times). If a reflection is needed, then for some $j$,
    $$f(x') = (x_j, x_{j - 1}, \dots, x_1, x_{2k}, \dots, x_{j + 1}).$$
    And so, $f(x)$ is in the same orbit as $f(x')$ under $D_{2k}$ (reflect then rotate $f(x)$ $j$ times).

    Now, let $x, x' \in X'_{2k}$ and suppose $f(x), f(x')$ are in the same orbit under $D_{2k}$. Let $f(x) = (x_1, x_2, \dots, x_{2k})$. Then, $f(x')$ may be obtained from $f(x)$ by a series of rotations and possibly a reflection. Similar to above, if no reflection is needed, then for some $i$,
    $$f(x') = (x_{2k-i}, x_{2k-i+1}, \dots, x_{2k}, x_1, \dots, x_{2k-i-1}).$$
    In this case, $x$ is in the same orbit as $x'$ under $D_n$ (rotate $x$ $(\sum_{m = 0}^i (1+ x_{2k - m}))$-many times). If a reflection is needed, then for some $j$,
    $$f(x') = (x_j, x_{j - 1}, \dots, x_1, x_{2k}, \dots, x_{j + 1}).$$
    This implies that $x$ is in the same orbit of $x'$ under $D_n$ (reflect and then rotate $x$ $(1 + \sum_{m = 1}^{j} (1+x_m))$-many times). Therefore, $x,x'$ are in the same orbit of $X_{2k}$ under $D_n$ if and only if $f(x), f(x')$ are in the same orbit of $C_{2k}(n-2k)$ under $D_{2k}$. Since each orbit of $X_{2k}$ under $D_n$ has some representative in  $X'_{2k}$, we can conclude that $|X_{2k}/D_{n}| = |C_{2k}(n-2k)/D_{2k}|$.

    For every $k$ in which $n-2k \ge 2k$, there exist $(2k)$-compositions of $n-2k$, and so, $|C_{2k}(n-2k)/D_{2k}|$ counts the number of distinguishable strings for such $k$. When $k=0$, there is one distinguishable cyclic string of $A$. Equality is thus proven for even $n$.
    
    If $n$ is odd, the same argument applies. For each $k>0$, let $X_{2k+1}$ be the set of cyclic length-$n$ strings of $A$ with $(2k+1)$-many $0$'s. Then,  we have $|X_{2k+1}/D_n| = |C_{2k+1}(n-2k-1)/D_{2k+1}|.$ For every $k$ in which $n-(2k+1) \geq 2k+1$, there exist $(2k+1)$-compositions of $n-2k-1$. When $k=0$, there is one distinguishable cyclic string.
\end{proof}

\subsection{Average size of MIS's}\label{avg size subsection}

Another area that has received attention (see~\cite{Barbosa}) is the size of maximal independent sets in grid-like graphs. As well as parameters such as the largest and smallest possible MIS of $G$, it is natural to consider the average size of an MIS. For a given graph $G$, we define the \textit{total MIS size} $$T(G):=\sum_{M \in \MIS(G)} |M|$$ 
and define the \textit{average MIS size} 
\[A(G) := \frac{T(G)}{|\MIS(G)|}.\]
We give explicit formulas for the total and average MIS sizes in some specific grid-like graphs.

For the $2 \times n$ grid graph $\Grid{2}{n}$, we can use the map $\psi$ from Definition~\ref{def:mappsi} to describe the behavior of $T(\Grid{2}{n})$. We set $X_n := \psi\left(\MIS(\Grid{2}{n})\right)$. By Lemma~\ref{map facts}, $X_n$ is the set of binary strings starting and ending with 1's and containing no adjacent 0's. Since every binary string corresponds to two MIS's, summing over all binary strings in $X_n$ yields $\frac{T(\Grid{2}{n})}{2}$.

\begin{lemma}\label{fibonacci convolution}
    We have $$T(\Grid{2}{n}) = \frac{2}{5}\left[nF(n+2) + (n+2)F(n)\right].$$
\end{lemma}
\begin{proof}
    Let $t(n)$ be the sum over all binary strings in $X_n$. We sum over all strings bit by bit. Let $i \in [n]$. When bit $i$ is a 1, the substring from bit 1 to bit $i$ (inclusive) is in $X_i$, and the substring from bit $i$ to bit $n$ (inclusive) is in $X_{n-i+1}$. Since every choice from $X_i$ and $X_{n-i+1}$ results in a binary string in $X_n$ with bit $i$ equal to 1, there are $|X_i||X_{n-i+1}|$ strings whose bit $i$ contributes to the sum. Thus, the sum over all binary strings in $X_n$ is $$\sum_{i=1}^n |X_i||X_{n-i+1}|.$$ 
    Applying Theorem~\ref{2 by n count} and Lemma~\ref{map facts},
    we have that $|X_k| = F(k)$, so $$t(n) = \sum_{i=1}^n F(i)F(n-i+1).$$ In particular, $t(n)$ is a convolution of Fibonacci numbers. These convolutions have been studied by various authors, with a formula for this convolution given in Equation 1.12 of \cite{Convolution}. Specifically, the authors use the convention $F(0) = 0$ and show that 
    \[5\sum_{i=0}^k F(i)F(k-i) = (k-1)F(k+1) + (k+1)F(k-1).\] Since $F(0) = F(k-k) = 0$, we can set $k = n+1$ in this formula and have \[\sum_{i=0}^{n+1} F(i) F(n-i+1) = \sum_{i=1}^{n} F(i) F(n-i+1).\] Thus we have
    \[\sum_{i=1}^n F(i)F(n-i+1)=\frac{1}{5}\left[nF(n+2)+(n+2)F(n)\right].\]
    Since $t(n)=\frac{T(\Grid{2}{n})}{2}$, the result follows.
\end{proof}

\gridaveragesize*

\begin{proof}
    By Theorem~\ref{2 by n count}, $|\MIS(\Grid{2}{n})| = 2F(n)$. So by Lemma~\ref{fibonacci convolution}, 
    \[A(\Grid{2}{n})=\frac{T(\Grid{2}{n})}{|\MIS(\Grid{2}{n})|}=\frac{1}{5}\left(\frac{nF(n+2) + (n+2)F(n)}{F(n)}\right) = \frac{1}{5}\left(\frac{nF(n+2)}{F(n)}+n+2\right).\] 
    Observing that $n \cdot \frac{F(n+2)}{F(n)}= n \varphi^2 + o(1)$, we have 
    \[A(\Grid{2}{n}) = \frac{1}{5}(n\varphi^2+n+2) + o(1) =\left(\frac{\varphi^2 + 1}{5}\right)n+\frac{2}{5} + o(1).\] 
    As 
    \[\frac{\varphi^2 + 1}{5} = \frac{\varphi}{\sqrt{5}} \approx 0.724,\] 
    the result follows.
\end{proof}

We note that the closed form in the limit is also a good approximation even for relatively small values of $n$: we computed that it is correct to 8 and 20 decimal places for $n=20$ and $n=50$, respectively.

We next consider the average size of MIS's in $\Band{2}{n}$. We use the map $\psi_c$ from Definition~\ref{def:mappsi} and set $Y_n := \psi_c(\MIS(\Band{2}{n}))$. As noted in Observation~\ref{cyclic map facts}, $Y_n$ is the set of cyclic length-$n$ binary strings with no adjacent 0's (recall that the first and last bit in a cyclic string are adjacent) and an even number of 1's. 

\bandaveragesize*

\begin{proof}
    We first count $|Y_n|$. Let $E_n$ be the set of non-cyclic length-$n$ binary strings with no adjacent 0's and an even number of 1's. Note that since the adjacency conditions of $E_n$ and $Y_n$ agree except with regard to the first and last bits, every string of $E_n$ is an element of $Y_n$ unless the string contains both the prefix 01 and the postfix 10. There are exactly $|E_{n-4}|$-many such strings, as for any string $s \in E_n \setminus Y_n$, the substring of $s$ strictly between indices 2 and $n-2$ is an element of $E_{n-4}$; note that two distinct strings $s,s' \in E_n \setminus Y_n$ must correspond to different strings in $E_{n-4}$, and each string in $E_{n-4}$ can be concatenated with prefix $01$ and postfix $10$ to form a string in $E_n \setminus Y_n$. Thus, $|Y_n| = |E_n| - |E_{n-4}|$. 

    To count $|E_n|$, observe that strings of $E_n$ begin with exactly one of the following prefixes: $11$, $101$, $011$, or $0101$. Because each prefix of length $i$ ends in a $1$ and has an even number of $1$'s, any element of $E_{n-i}$ could be appended after the prefix to obtain a string of $E_n$. Thus, $|E_n| = |E_{n-2}| + 2|E_{n-3}| + |E_{n-4}|$ with $|E_1| = 1, |E_2| = 1, |E_3| = 3, |E_4| = 4$. Appearing as A093040 in Sloane's On-line Encyclopedia of Integer Sequences \cite{OEIS_A093040}, several closed-form formulas for $|E_n|$ have been observed. We will use the following form, contributed by Detlefs (and quickly verifiable by induction): for all $n\geq 1$, $$|E_n| = \left\lfloor \frac{F(n+3)}{2} \right\rfloor - \left\lfloor \frac{F(n+1)}{2} \right\rfloor.$$

    Let $t_c(n)$ be the sum over all binary strings in $Y_n$. Similar to the methodology in Lemma \ref{fibonacci convolution}, we sum over all strings bit by bit. The number of strings of $Y_n$ with a $1$ at index $i$ is the number of ways to form a length-$n$ binary cyclic string which satisfies the conditions on strings in $Y_n$ and contains a $1$ at index $i$. Note that every such string $s$ corresponds to a non-cyclic binary string $s'$ of length $n-1$ with an odd number of $1$'s and no adjacent $0$'s. (Specifically, we form $s'$ from the bits of $s$, written in the order $i+1, i+2, \dots,n,1,2 \dots, i-1$.) Thus, the number of elements of $Y_n$ with a $1$ in bit $i$ is equal to the number of length-$(n-1)$ non-cyclic binary strings satisfying the above conditions.

    Observe that $F(n+2)$ counts the size of the set $B_n$ of $n$-bit non-cyclic binary strings with no repeating 0's. 
    Indeed, we have $|B_1| = 2 = F(3)$ and $|B_2| = 3 = F(4)$; for $n \geq 3$, note that every string in $B_n$ is either obtained by appending a string in $B_{n-1}$ to the prefix $1$ or a string in $B_{n-2}$ to the prefix $01$, so $|B_n| = |B_{n-1}| + |B_{n-2}|$ obeys the Fibonacci recurrence with the correct initial conditions. Thus, $F(n+1)-|E_{n-1}|$ is the number of length-$(n-1)$ non-cyclic binary strings with no repeating $0$'s and an odd number of $1$'s. Note that the above count does not depend upon the value of index $i$: because the binary strings in $Y_n$ are cyclic, the number of ways to fill in the string around the $1$ in index $i$ is the same for all $i$. Thus, $t_c(n) = n(F(n+1)-|E_{n-1}|)$, and \[A(\Band{2}{n}) = \frac{2t_c(n)}{2|Y_n|} = \frac{n(F(n+1)-|E_{n-1}|)}{|E_n| - |E_{n-4}|} = n \cdot \frac{F(n+1) - \left\lfloor \frac{F(n+2)}{2} \right\rfloor + \left\lfloor \frac{F(n)}{2} \right\rfloor}{\left\lfloor \frac{F(n+3)}{2} \right\rfloor - \left\lfloor \frac{F(n+1)}{2} \right\rfloor - \left\lfloor \frac{F(n-1)}{2} \right\rfloor + \left\lfloor \frac{F(n-3)}{2} \right\rfloor}.\]

    We may rewrite \[\frac{F(n+1) - \left\lfloor \frac{F(n+2)}{2} \right\rfloor + \left\lfloor \frac{F(n)}{2} \right\rfloor}{\left\lfloor \frac{F(n+3)}{2} \right\rfloor - \left\lfloor \frac{F(n+1)}{2} \right\rfloor - \left\lfloor \frac{F(n-1)}{2} \right\rfloor + \left\lfloor \frac{F(n-3)}{2} \right\rfloor} = \frac{F(n+1) -  \frac{F(n+2)}{2} +  \frac{F(n)}{2} + O(1)}{ \frac{F(n+3)}{2}  -  \frac{F(n+1)}{2}  - \frac{F(n-1)}{2} +\frac{F(n-3)}{2} + O(1)}\] \[= \frac{\varphi^{n+1} -  \frac{\varphi^{n+2}}{2} +  \frac{\varphi^{n}}{2} + O(1)}{ \frac{\varphi^{n+3}}{2}  -  \frac{\varphi^{n+1}}{2}  - \frac{\varphi^{n-1}}{2} +\frac{\varphi^{n-3}}{2} + O(1)}\]
    since $F(n) = \frac{1}{\sqrt{5}}\varphi^n + o(1)$. We now simplify a bit further, multiplying numerator and denominator by $\frac{2}{\varphi^{n-3}}$:
    \[\frac{\varphi^{n+1} -  \frac{\varphi^{n+2}}{2} +  \frac{\varphi^{n}}{2} + O(1)}{ \frac{\varphi^{n+3}}{2}  -  \frac{\varphi^{n+1}}{2}  - \frac{\varphi^{n-1}}{2} +\frac{\varphi^{n-3}}{2} + O(1)} = \frac{2\varphi^{4} -  \varphi^{5} +  \varphi^{3} + O(\varphi^{-n+3})}{ \varphi^{6}  -  \varphi^{4}  - \varphi^{2} + 1 + O(\varphi^{-n+3})}.\]
    Thus,
    \[A(\Band{2}{n}) = n \cdot \frac{2\varphi^{4} -  \varphi^{5} +  \varphi^{3} + O(\varphi^{-n+3})}{ \varphi^{6}  -  \varphi^{4}  - \varphi^{2} + 1 + O(\varphi^{-n+3})} = \frac{n \left(2\varphi^{4} -  \varphi^{5} +  \varphi^{3}\right)}{ \varphi^{6}  -  \varphi^{4}  - \varphi^{2} + 1} + o(1) = \left(\frac{\varphi}{\sqrt{5}}\right)n + o(1).\]
    It is straightforward but tedious to verify the last equality by hand using identities on powers of $\varphi$; we obtained the reduced form computationally.
\end{proof}

Thus, we have from Theorem~\ref{Grid Average Size Thm} and Theorem~\ref{Band Average Size Thm} that \[\lim_{n \to \infty} \left[ A(\Grid{2}{n})-A(\Band{2}{n}) \right]=\frac{2}{5}.\]
In particular, the discrepancy between $A(\Grid{2}{n})$ and $A(\Band{2}{n})$ does not vanish as $n$ grows (although, for large $n$, this difference is very small in comparison to the total values of $A(\Grid{2}{n})$ and $A(\Band{2}{n})$).

We also comment that Barbosa and Cappelle~\cite{Barbosa} showed that the minimum size of an MIS in $\Band{2}{n}$ is $2 \lceil \frac{n}{4} \rceil$ and that the maximum size of an MIS is $2 \lfloor \frac{n}{2} \rfloor$. Since $\frac{\varphi}{\sqrt{5}} \approx 0.724$, Theorem~\ref{Band Average Size Thm} shows that in the limit, the average size of an MIS is slightly closer to the minimum possible size than to the maximum possible size. In the following theorem, we expand on the description of MIS's in $\Band{2}{n}$ and enumerate MIS's of each fixed size in $\Band{2}{n}$ and $\Mobius{2}{n}$. Analogous to the ideas of Observation~\ref{cyclic map facts}, the image of $M \in \MIS(\Mobius{2}{n})$ under $\psi_c$ must be a cyclic string with no adjacent zeros and an odd total number of $1$'s.

\begin{thm}\label{thm:B2xn sizes}
    Let $n \geq 3$ and $0 < r \le n$. Then, there are \[2\left(\binom{r}{n-r}+\binom{r-1}{n-r-1}\right)\] $r$-vertex elements of $\MIS(\Band{2}{n})$ if $r$ is even, or of $\MIS(\Mobius{2}{n})$ if $r$ is odd. 
\end{thm}

\begin{proof}
    First suppose $r$ is even, and consider the set $R$ of elements of $\MIS(\Band{2}{n})$ with exactly $r$ vertices. Note that elements of $R$ can be paired by reflecting an MIS in $R$ over the horizontal axis between rows 1 and 2. Form $R'\subset R$ by deleting exactly one element from each pair in $R$. Note that $|R| = 2|R'|$ and recall the function $\psi_c$ from Definition~\ref{def:mappsi}. Then, $\psi_c(R')$ is precisely the set of elements in $Y_n$ containing exactly $r$-many $1$'s. Consider a string $s \in \psi_c(R')$; then, there are $r$ many 1's and $(n-r)$-many 0's in $s$. Because $s$ contains no two adjacent $0$'s, there must be at least one 1 separating every two closest $0$'s, and $s(1)$ and $s(n)$ cannot both be $0$. We build $s$ by first placing the $r$-many $1$'s and then placing the $(n-r)$-many $0$'s. If $s(1)=1$, we can choose any of the $r$-many $1$'s to place a $0$ to the right of. Otherwise, $s(1)=0$ and $s(2)=1=s(n)$, and we can choose any of the $(r-1)$-many $1$'s after the $01$-prefix of $s$ to place a $0$ to the left of. Thus, $|R'|=\binom{r}{n-r}+\binom{r-1}{n-r-1}$ as desired.
    
    If $r$ is odd, essentially identical reasoning can be used to show the identity holds for $\MIS(\Mobius{2}{n})$. We note that $\Band{2}{n}$ has no MIS of odd size, $\Mobius{2}{n}$ has no MIS of even size, and (as the identity implies) neither $\Band{2}{n}$ nor $\Mobius{2}{n}$ have any $r$-vertex MIS's if $n > 2r$.  
\end{proof}

Theorem~\ref{thm:B2xn sizes} illustrates strongly similar behavior in the distribution of sizes of MIS in $\Band{2}{n}$ and $\Mobius{2}{n}$. This may suggest that $A(\Mobius{2}{n})$ should have a similar value (and be understandable in similar fashion) to $A(\Band{2}{n})$. This is indeed the case.

\mobiusaveragesize*

\begin{proof}
Let $Z_n := \psi_c(\MIS(\Mobius{2}{n}))$. Note that $Z_n$ is the set of cyclic length-$n$ binary strings with no adjacent $0$'s and an odd number of $1$'s. As in Theorem~\ref{Band Average Size Thm}, we let $E_n$ be the set of non-cyclic length-$n$ binary strings with no adjacent $0$'s and an even number of $1$'s. We define $O_n$ to be the set of non-cyclic length-$n$ binary strings with no adjacent $0$'s and an odd number of $1$'s. Thus, $|E_n| + |O_n|$ is equal to the total number of length-$n$ binary strings with no adjacent $0$'s; we recall that this total number is $F(n+2)$, and that $|E_n| = \left\lfloor \frac{F(n+3)}{2} \right\rfloor - \left\lfloor \frac{F(n+1)}{2} \right\rfloor$. Thus, $|O_n| = F(n+2) - \left\lfloor \frac{F(n+3)}{2} \right\rfloor + \left\lfloor \frac{F(n+1)}{2} \right\rfloor$. As in the proof of Theorem~\ref{Band Average Size Thm}, the only elements of $O_n$ which are not in $Z_n$ are those strings beginning with prefix $01$ and ending with postfix $10$, which implies $|Z_n| = |O_n| - |O_{n-4}|$. 

Let $t_c(n)$ be the sum over all binary strings in $Z_n$. Note that $2 t_c(n) = T(\Mobius{2}{n})$ and $2 |Z_n| = |\MIS(\Mobius{2}{n})|$, so 
\[A(\Mobius{2}{n}) = \frac{t_c(n)}{|Z_n|}.\]
As in Lemma~\ref{fibonacci convolution} and Theorem~\ref{Band Average Size Thm}, we determine $t_c(n)$ by summing over strings bit-by-bit. The number of strings in $Z_n$ with a $1$ at index $i$ is the number of ways to complete a cyclic string around the $1$ at index $i$. This is the same as the number of non-cyclic strings of length $n-1$ with an even number of $1$'s and no repeating $0$'s. This analysis holds for each of the $n$ bits, so $t_c(n) = n |E_{n-1}|$. Applying the formulas for $|Z_n|$, $|E_n|$, and $|O_n|$, we thus have
\[A(\Mobius{2}{n}) = n \cdot \frac{\left\lfloor \frac{F(n+2)}{2} \right\rfloor - \left\lfloor \frac{F(n)}{2} \right\rfloor}{F(n+2) - \left\lfloor \frac{F(n+3)}{2} \right\rfloor + \left\lfloor \frac{F(n+1)}{2} \right\rfloor - F(n-2) + \left\lfloor \frac{F(n-1)}{2} \right\rfloor - \left\lfloor \frac{F(n-3)}{2} \right\rfloor}.\]
We estimate this value using that $F(n) = \frac{1}{\sqrt{5}}\varphi^n + o(1)$ and $\left\lfloor \frac{F(n)}{2}\right\rfloor = \frac{1}{\sqrt{5}}\frac{\varphi^n}{2} + O(1)$, so
\[A(\Mobius{2}{n}) = n \cdot \frac{ \frac{\varphi^{n+2}}{2}  -  \frac{\varphi^{n}}{2} + O(1)}{\varphi^{n+2} -  \frac{\varphi^{n+3}}{2}  + \frac{\varphi^{n+1}}{2}  - \varphi^{n-2} +  \frac{\varphi^{n-1}}{2}  -  \frac{\varphi^{n-3}}{2} + O(1)}.\]
Multiplying numerator and denominator by $\frac{2}{\varphi^{n-3}}$, we obtain

\[A(\Mobius{2}{n}) = n \cdot \frac{\varphi^{5}  -  \varphi^{3}  + O\left(\frac{1}{\varphi^{n-3}}\right)}{2\varphi^{5}  -  \varphi^{6}   + \varphi^{4}   - 2\varphi +  \varphi^{2}   -  1 + O\left(\frac{1}{\varphi^{n-3}}\right)}\]
\[= \frac{n\left(\varphi^{5}  -  \varphi^{3}\right)}{2\varphi^{5}  -  \varphi^{6}   + \varphi^{4}   - 2\varphi +  \varphi^{2}   -  1} + o(1) = \left(\frac{\varphi}{\sqrt{5}}\right)n + o(1).\]
As in Theorem~\ref{Band Average Size Thm}, the last inequality is non-obvious, but follows either from by-hand verification using identities on powers of $\varphi$, or by computer calculation.
\end{proof}

\section*{Acknowledgements}

The results in this paper are based on work completed at the Math REU at Iowa State University in 2024. The authors would like to thank Steve Butler and Kate Lorenzen for their work in organizing this program and their mentorship contributions throughout.

\bibliographystyle{plain} 
\bibliography{refs} 
\nocite{}
\end{document}